\documentclass[12pt]{amsart}
\usepackage{amsmath,amssymb}

\newtheorem{theorem}{Theorem}

\renewcommand{\L}{\mathcal{L}}
\def\cset#1#2{\{#1 : #2\}}
\def\seq#1#2{\langle#1 : #2\rangle}
\def\tp{\mathrm{tp}}
\def\Aut{\mathrm{Aut}}

\begin{document}
\title{The strict order property and generic automorphisms}
\author{Hirotaka Kikyo}
\address{\hskip-\parindent
	Hirotaka Kikyo\\
	Department of Mathematical Sciences\\
	Tokai University\\
	1117 Kitakaname, Hiratsuka, 259-1292\\
	Japan}
\email{kikyo@ss.u-tokai.ac.jp}

\author{Saharon Shelah*}
\thanks{%
* The second author is supported by the United States-Israel Binational
Science Foundation. Publication 748.
}
\address{\hskip-\parindent
         Saharon Shelah\\
         Institute of Mathematics\\
         The Hebrew University\\ 
         Jerusalem 91904, Israel}
\email{shelah@sundial.ma.huji.ac.il}

\begin{abstract}
If $T$ is an model complete theory with the strict order property, 
then the theory of the models of $T$ with an automorphism
has no model companion.
\end{abstract}
\maketitle

\section{Introduction}

Given a model complete theory $T$ in a language $\L$, we
consider the (incomplete) theory $T_\sigma =
T \cup \{\text{``$\sigma$ is an $\L$-automorphism''}\}$ in the
language $\L_\sigma = \L \cup \{\sigma\}$.
For $M$ a model of $T$, and $\sigma \in \Aut_\L(M)$ we call
$\sigma$ a generic automorphism of $M$ if $(M,\sigma)$ is 
an existentially closed model of $T_{\sigma}$.
A general problem is to find necessary and sufficient conditions 
on $T$ for the class of existentially models of $T_{\sigma}$ 
to be elementary,
namely to be the class of models of some first
order theory in $\L_\sigma$. This first order theory, if it exists, is
denoted $TA$, and it is the model companion of $T_\sigma$. 
This problem seems to be a difficult problem
even if we assume $T$ to be stable 
\cite{Chatzidakis}, \cite{KP}, \cite{Macintyre}.
Generic automorphisms in the sense of this paper were first studied by
Lascar \cite{Lascar2}. 
The work of Chatzidakis and
Hrushovski \cite{Chatzidakis-Hrushovski} on the case where $T$ is the theory
$ACF$ of algebraically closed fields renewed interest in the topic
and Chatzidakis and Pillay studied general properties of $TA$ 
for stable $T$ \cite{ChP}.

Kudaibergenov proved that if $TA$ exists then $T$ eliminates the
quantifier ``there exists infinitely many''.
Therefore, if $T$ is stable and $TA$ exists then
$T$ does not have the fcp.
Pillay conjectured that if $T$ has the fcp then $TA$ does not
exist after the first author observed that the theory of
random graphs does not have $TA$. 
The first author then proved that if $TA$ exists and $T$ does not 
have the independence property then $T$ is stable,
and if $TA$ exists and $T_\sigma$ has the amalgamation property
then $T$ is stable \cite{K}. 
The latter fact covers the case of the random graphs.
The present paper extends the former case.

So the theorem here shows that model complete theories with $T_\sigma$
having a model companion are ``low'' in the hierarchy of
\emph{classification theory}; 
previous results have shown it cannot be in some
intermediate positions.

For other examples,
Hrushovski observed that there are no $TA$ for
$ACFA$ and for the theory of pseudo-finite fields $Psf$
(unpublished).
His argument depends heavily on field theory.
$ACFA_\sigma$ does not have the amalgamation property
but it is not known if $Psf_\sigma$ has the amalgamation property.

In the rest of the paper, 
small letters $a$, $b$, $c$, etc.\ denote finite tuples
and $x$, $y$, etc.\ denote finite tuples of variables.
If $a$ is a tuple of elements and $A$ a set of elements,
$a \in A$ means that each element of $a$ belongs to $A$.

\section{Main Theorem}

\begin{theorem} Let $T$ be a model complete theory in a language $\L$
and $\sigma$ a new unary function symbol.
If $T$ has a model whose theory has the strict order property then
$T_\sigma$ has no model companion.
\end{theorem}

\begin{proof}
Let $M_0$ be a model of $T$ with the strict order property.
So, there are $\L$-definable partial order $<$ on $k$-tuples in $M_0$
for some $k$ and 
a sequence $\seq{a_i}{i<\omega}$ of $k$-tuples in $M_0$ such that 
$a_i < a_j$ for $i < j < \omega$.
By Ramsey's Theorem, we can assume that
$\seq{a_i}{i<\omega}$ is an $\L$-indiscernible sequence in $M_0$.
Also, we can assume that there is an $\L$-automorphism $\sigma_0$ of $M_0$
such that $\sigma_0(a_i) = a_{i+1}$.
So, $(M_0, \sigma_0)$ is a model of $T_\sigma$.

Now by way of contradiction,
suppose that $T_\sigma$ has a model companion, say $TA$.
Extend $(M_0, \sigma_0)$ to a model $(N,\sigma)$ of $TA$.
$N$ is an $\L$-elementary extension of $M_0$ since $T$ is
model complete.
We can assume that $(N, \sigma)$ is sufficiently saturated.
In the rest of the proof, we work in $(N, \sigma)$.

Consider the partial type $p(x) = \cset{a_i < x}{i<\omega}$ and
let $\psi(x) \equiv 
\exists y (a_0 < \sigma(y) \land \sigma(y) < y \land y < x)$.

\medskip

\noindent
{\bf Claim.} In $(N, \sigma)$,
\begin{enumerate}
\renewcommand{\labelenumi}{(\arabic{enumi})}
\item $p(x) \vdash \psi(x)$, and
\item if $q(x)$ is a finite subset of $p(x)$ then 
$q(x) \not\mathop{\vdash} \psi(x)$.
\end{enumerate}

\medskip

If this claim holds, then it contradicts the saturation of $(N, \sigma)$.

We first show (2).
Let $n^*$ be such that $q(x) \subset \cset{a_i< x}{i < n^*}$.
Then $a_{n^*}$ satisfies $q(x)$.
Suppose $a_{n^*}$ satisfies $\psi(x)$. Let $b \in N$ be such that
$a_0 < \sigma(b) < b < a_{n^*}$.
By $a_0 < \sigma(b)$, 
we have $a_{n^*} = \sigma^{n^*}(a_0) < \sigma^{n^*+1}(b)$.
By $\sigma(b) < b < a_{n^*}$, we have
\[
\sigma^{n^*+1}(b) < \sigma^{n^*}(b) < \cdots < \sigma(b) < b < a_{n^*}.
\]
By transitivity, we get $a_{n^*} < a_{n^*}$, which is a contradiction.

Now we turn to a proof of (1).
Suppose $c \in p(N)$.
Let $M$ be such that 
$a_0, c \in M$, $|M| = |T|$, and
$(M,\sigma|M)$ is an $\L(\sigma)$-elementary
substructure of $(N,\sigma)$.

For each $d \in p(N)$, 
let $\Psi(d)$ be the set of $\L(M)$-formulas $\varphi(x)$ satisfied 
in $N$ by some tuple $d'$ such that $d'\in p(N)$ and $d' < d$.
Here, $\L(M)$-formulas are the formulas in $\L$ with parameters in $M$.

Note that
if $d_1, d_2 \in p(N)$ and $d_2 < d_1$, then $\Psi(d_2) \subseteq \Psi(d_1)$,
and by compactness,
if $d_1, d_2 \in p(N)$ then there is $d_3 \in p(N)$ such that
$d_3 < d_1$ and $d_3 < d_2$.

Since $N$ is sufficiently saturated,
there is $c^* \in p(N)$ such that
if $d \in p(N)$ and $d < c^*$ then $\Psi(d) = \Psi(c^*)$.
We can also assume that $c^* < c$.
Since the sets $p(N)$ and $M$ are invariant under $\sigma$,
$\Psi(c^*)$ is also invariant under $\sigma$, which means,
for any $\L$-formula $\varphi(x,y)$ and a tuple $a \in M$,
$\varphi(x,a) \in \Psi(c^*)$ if and only if
$\varphi(x,\sigma(a)) \in \Psi(c^*)$.

Now choose $b_1 \in p(N)$ such that $b_1 < c^*$ and 
consider $q_1(x) = \tp_{\L}(b_1/M)$. Then $q_1(x) \subseteq \Psi(c^*)$.
Let $\sigma(q_1(x))$ be the set of formulas $\varphi(x,\sigma(a))$
such that $\varphi(x,a) \in q_1(x)$ where $\varphi(x,y)$ is a 
formula in $\L$ and $a \in M$.
Since $\Psi(c^*)$ is invariant under $\sigma$, we have
$\sigma(q_1(x)) \subseteq \Psi(c^*)$.
By the choice of $c^*$, $\Psi(c^*) = \Psi(b_1)$ and thus
$\sigma(q_1(x)) \subseteq \Psi(b_1)$.
By the definition of $\Psi(b_1)$ and by compactness,
there is $b_2 \in p(N)$ such that $b_2 < b_1$ and
$b_2$ realizes $\sigma(q_1(x))$.

Since $\sigma(q_1(x))$ is a complete $\L$-type over $M$,
there are $\L$-substructure $M'$ of $N$
and an $\L$-automor\-phism $\tau$ of $M'$ such that
$Mb_1b_2 \subset M'$, $\tau(b_1) = b_2$ and $\tau|M = \sigma|M$.
Now we have,
\[
(M',\tau) \models a_0 < \tau(b_1) < b_1 < c
\]
Since $(M,\sigma|M)$ is a model of $TA$, it is an existentially closed
model of $T_\sigma$.
Note that the partial order $<$ is definable by an existential $\L$-formula
modulo $T$.
So, the formula $a_0 < \sigma(y) < y < c$ has a solution in $(M,\sigma|M)$.
Hence, we have $(M,\sigma|M) \models \psi(c)$.
This proves Claim (1) and we are done.
\end{proof}

\end{document}